\documentclass[12pt,draft]{amsart}

\makeatletter

\theoremstyle{plain}    
\newtheorem{thm}{Theorem}[section]
\numberwithin{equation}{section} 
\numberwithin{figure}{section} 
\theoremstyle{plain}    
\newtheorem*{thm*}{Theorem} 
\theoremstyle{plain}    
\newtheorem{cor}[thm]{Corollary} 
\theoremstyle{plain}    
\newtheorem{lem}[thm]{Lemma} 
\theoremstyle{plain}    
\newtheorem{prop}[thm]{Proposition} 
\theoremstyle{definition}
\newtheorem{defn}[thm]{Definition}
\theoremstyle{remark}

\theoremstyle{remark}

\theoremstyle{remark}    

\theoremstyle{remark}    

\theoremstyle{definition}  

\theoremstyle{remark}
  \newtheorem*{acknowledgement*}{Acknowledgement} 

\theoremstyle{plain}    

\theoremstyle{plain}    
\theoremstyle{plain}    
\theoremstyle{plain}    
\theoremstyle{definition}

\theoremstyle{remark}

\theoremstyle{remark}    

\theoremstyle{remark}    

\theoremstyle{plain}    

\newcommand{\TAim}{\operatorname{T(A)_{IM}}}

\usepackage{amssymb} \usepackage{amscd}

\makeatother

\begin{document}

\title[Connes' Embedding Problem]{Connes' Embedding Problem and Lance's WEP}

\author{Nathanial P. Brown}

\address{Department of Mathematics, Penn State University, State
College, PA 16802}

\email{nbrown@math.psu.edu}

\thanks{Supported by an NSF Postdoctoral
Fellowship.  2000 MSC number: 46L05. }

\begin{abstract}
A II$_1$-factor embeds into the ultraproduct of the hyperfinite
II$_1$-factor if and only if it satisfies the von Neumann algebraic
analogue of Lance's weak expectation property (WEP).  This note gives
a self contained proof of this fact.
\end{abstract}

\maketitle

\section{Introduction}

On page 105 in \cite{connes:classification} Connes suggested that
every separable II$_1$-factor ought to be embeddable into the
ultraproduct, $R^{\omega}$, of the hyperfinite II$_1$-factor $R$.
Largely due to work of Kirchberg, Voiculescu and, most recently,
Haagerup this seemingly technical question has received more and more
attention in recent years.  Indeed, Kirchberg proved in 
\cite{kirchberg:invent} that this
problem can be reformulated in an unexpected variety of ways (see
\cite{ozawa:QWEP} for a wonderful exposition of Kirchberg's work),
this problem turns out to be a necessary condition for Voiculescu's
`Unification Problem' (i.e.\ if the microstates and non-microstates
approaches to free entropy yield the same quantity then every
II$_1$-factor is embeddable) and, finally, Haagerup has shown that
this problem is nearly sufficient for resolving the relative invariant
subspace problem for II$_1$-factors (he showed that every operator in
an embeddable II$_1$-factor which satisfies a mild non-degeneracy
condition has invariant subspaces -- see
\cite{haagerup:invariantsubspaces}).

In \cite{lance:nuclear} Lance introduced the {\em weak expectation
property} (WEP) for C$^*$-algebras.  Blackadar shifted the point of view to
von Neumann algebras with the following definition.

\begin{defn}
Let $M \subset B(H)$ be a von Neumann algebra acting on some Hilbert
space $H$ and let $A \subset M$ be a weakly dense C$^*$-subalgebra.
Then $M$ has a {\em weak expectation relative to $A$} if there exists
a unital, completely positive map $\Phi : B(H) \to M$ such that
$\Phi(a) = a$ for all $a \in A$.
\end{defn}

This notion was inspired by injectivity; $M \subset B(H)$ is {\em
injective} if there exists a unital, completely positive map $\Phi :
B(H) \to M$ such that $\Phi(x) = x$ for all $x \in M$.

It follows from Arveson's Extension Theorem that a C$^*$-algebra $A$
has the WEP if and only if the enveloping von Neumann algebra $A^{**}$
has a weak expectation relative to $A \subset A^{**}$.  

In \cite{brown:invariantmeans} we observed that the W$^*$-version of
the WEP is closely related to  Connes' embedding problem. 

\begin{thm}
\label{thm:mainthm}
Let $M$ be a separable II$_1$-factor.  Then $M$ is embeddable into
$R^{\omega}$ if and only if $M$ has a weak expectation relative to
some weakly dense subalgebra.
\end{thm}

A simple corollary of this result states that many well known
II$_1$-factors which are ``far from being hyperfinite'' (in the sense
that they exhibit vastly different properties than $R$ -- no Cartan
subalgebras, prime,  property T, etc.) are
in fact built out of $R$ in a way which naturally mixes von Neumann
algebraic and operator space notions.  More precisely, we have the
following approximation property.

\begin{cor}
Let $M \subset R^{\omega}$ be a II$_1$-factor, ${\mathfrak F} \subset
M$ be a finite set and $\epsilon > 0$ be given.  Then there exists a
subspace $X \subset M$ such that $X$ nearly contains ${\mathfrak F}$
(within $\epsilon$ in 2-norm) and $X \cong R$ (as operator systems).
\end{cor}

In other words, free group factors (and $L(\Gamma)$ for any other
residually finite group) are the weak closure of (operator space
isomorphic) copies of $R$.

The purpose of this note is to give  self contained proofs of these
results as some details do not appear in
\cite{brown:invariantmeans}.  The proof turns out to be fairly
elementary but relies on a mixture of classical ideas (invariant
means) some new aspects of (finite) representation theory of
C$^*$-algebras and a bit of trickery. 

Throughout this paper $A$ will denote a separable unital
C$^*$-algebra.  Separability is really not necessary, but it is
convenient. We will use the abbreviation u.c.p.\ for unital completely
positive maps.  If $\tau$ is a state on a C$^*$-algebra $A$ then
$\pi_{\tau} : A \to B(L^2(A,\tau))$ will denote the GNS
representation.  Note that if $M$ is a II$_1$-factor in standard form
(i.e.\ acting, via GNS, on the $L^2$-space coming from its unique
trace) and $\pi : A \to M \subset B(L^2(M))$ is a $*$-homomorphism
with weakly dense range then we may, thanks to uniqueness of GNS
representations, identify $\pi$ with the GNS representation of $A$
coming from $\tau \circ \pi$, where $\tau$ is the unique trace on $M$.

Finally, recall that if $R$ denotes the hyperfinite II$_1$-factor and
$\omega \in \beta({\mathbb N})\backslash {\mathbb N}$ is a free
ultrafilter then the ultraproduct $R^{\omega}$ is defined to be
$l^{\infty}(R) = \{ (x_n) : x_n \in R, \sup_n \| x_n \| < \infty \}$
modulo the ideal $I_{\omega} = \{ (x_n) : \lim_{n \to \omega} \| x_n
\|_2 = 0 \}$, where $\| x \|_2^2 = \tau(x^* x)$ and $\tau$ is the
unique trace on $R$.  It turns out that $R^{\omega}$ is a
II$_1$-factor with tracial state $\tau_{\omega} ((x_n)) = \lim_{n \to
\omega} \tau(x_n)$.

\section{Invariant Means on C$^*$-algebras} 

In \cite[Remark 5.35]{connes:classification} Connes points out that a
hypertrace can be regarded as the analogue of an invariant mean on a
group.  We essentially take this as the definition of an invariant
mean on a C$^*$-algebra. 

\begin{defn} Let $A \subset B(H)$ be a C$^*$-algebra.  A (tracial)
state $\tau$ on $A$ is called an {\em invariant mean} if there exists
a state $\psi$ on $B(H)$ which is (1) invariant under the action of
the unitary group of $A$ on $B(H)$ (i.e.\ $\psi(uTu^*) = \psi(T)$ for
all $T \in B(H)$ and unitaries $u \in A$) and (2) extends $\tau$
(i.e.\ $\psi|_A = \tau$).  We will denote by $\TAim$ the set of all
invariant means on $A$.
\end{defn}

The main result of this section gives an important characterization of
invariant means.  There are several other ways to characterize
invariant means (cf.\ \cite[Theorem 3.1]{brown:invariantmeans},
\cite[Theorem 6.1]{ozawa:QWEP}) but we only present the ones we need.
The main step in the proof ((1) $\Longrightarrow$ (2)) is essentially
due to Connes in the unique trace case and Kirchberg in general.  We
will isolate the main technical aspects in a lemma.

Below, Tr$(\cdot)$ will denote the canonical (unbounded) trace on
$B(H)$ and, if $H$ is finite dimensional, tr$(\cdot)$ will denote the
(unique) tracial state on $B(H)$.  Also, ${\mathcal T} \subset B(H)$
will be the trace class operators (i.e.\ the predual of $B(H)$) and
$\| \cdot \|_{1,Tr}$ (resp.\ $\| \cdot \|_{2,Tr}$) will denote the
L$^1$-norm (resp.\ L$^2$-norm) on ${\mathcal T}$.  Recall that the
Powers-St{\o}rmer inequality states that if $h,k \in {\mathcal T}$ are
positive then $\| h - k \|_{2,Tr}^2 \leq \| h^2 - k^2 \|_{1,Tr}$.  In
particular, if $u \in B(H)$ is a unitary and $h \geq 0$ has finite
rank then $\| uh^{1/2} - h^{1/2}u \|_{2,Tr} = \| uh^{1/2}u^* - h^{1/2}
\|_{2,Tr} \leq \| uhu^* - h \|_{1,Tr}^{1/2}$.

\begin{lem}
\label{thm:ozawalemma} Let $h \in B(H)$ be a positive, finite rank 
operator with rational eigenvalues and Tr$(h) = 1$.  Then there exists
a u.c.p.\ map $\phi : B(H) \to M_q({\mathbb C})$ such that
tr($\phi(T)) = {\rm Tr}(hT)$ for all $T \in B(H)$ and
$|$tr$(\phi(uu^*) - \phi(u)\phi(u^*))| < 2\| uhu^* - h \|_1^{1/2}$ for
every unitary operator $u \in B(H)$.
\end{lem}

\begin{proof} This proof is taken directly from the proof
of \cite[Theorem 6.1]{ozawa:QWEP} which, in turn, is based on work of
Haagerup.

Let $v_1, \ldots, v_k \in H$ be the eigenvectors of $h$ and
$\frac{p_1}{q}, \ldots, \frac{p_k}{q}$ the corresponding eigenvalues.
Note that $\sum p_j = q$ since Tr$(h) = 1$. Let $\{ w_m \}$ be any
orthonormal basis of $H$ and consider the following orthonormal subset
of $H\otimes H$: $$\{ v_1 \otimes w_1, \ldots, v_1 \otimes w_{p_1} \}
\cup \{ v_2 \otimes w_1, \ldots, v_2 \otimes w_{p_2} \} \cup \ldots 
\cup \{ v_k \otimes w_1, \ldots, v_k \otimes w_{p_k} \}.$$ Let 
$P \in B(H\otimes H)$ be the orthogonal projection onto the span of
these vectors.  We encourage the reader to write down the matrix of $P
(T\otimes 1) P$ (in the basis above), for an arbitrary $T \in B(H)$.
Indeed, having done so the following facts become easy to verify.

\begin{enumerate}
\item $\frac{1}{q} Tr(P(T\otimes 1)P) = \frac{1}{q}\big(\sum_{i =
1}^{k} p_i <Tv_i, v_i > \big) = \sum_{i = 1}^{k} <Thv_i, v_i> =
Tr(Th)$.

\item $\frac{1}{q} Tr(P(T\otimes 1)P(T^* \otimes 1)P) = \sum_{i,j =
1}^{k} |<Tv_j, v_i >|^2 \min\{p_i, p_j\}$. 
\end{enumerate}

Now we encourage the reader to write down the matrices of $h^{1/2}T$,
$h^{1/2}T^*$ and $h^{1/2}Th^{1/2}T^*$ (in any orthonormal basis which
begins with $\{ v_1, \ldots, v_k \}$).  Having done so one immediately
sees that, letting $T_{i,j} = <Tv_j, v_i>$, $$Tr(h^{1/2}Th^{1/2}T^*) =
\sum_{i,j = 1}^{k} \frac{1}{q}(p_i p_j)^{1/2} |T_{i,j}|^2.$$ 
Hence, if we define a u.c.p.\ map $\phi : B(H) \to M_q({\mathbb C})$
by $\phi(T) = P(T\otimes 1)P$ then $tr(\phi(T)) = Tr(hT)$ for all $T
\in B(H)$ and, moreover, we have the following estimates:
\begin{eqnarray*}
|Tr(h^{1/2}Th^{1/2}T^*) - tr(\phi(T)\phi(T^*))|  & = & \sum_{i,j = 1}^{k}
\frac{1}{q}|T_{i,j}|^2 \bigg((p_i p_j)^{1/2} - \min\{p_i, p_j\} \bigg)\\ 
& \leq &  \sum_{i,j = 1}^{k} \frac{1}{q}|T_{i,j}|^2 p_i^{1/2}|p_i^{1/2} 
          -  p_j^{1/2}|\\
& \leq & \bigg(\sum_{i,j = 1}^{k} \frac{1}{q}|T_{i,j}|^2 p_i   \bigg)^{1/2}  
         \bigg(\sum_{i,j = 1}^{k} \frac{1}{q}|T_{i,j}|^2 (p_i^{1/2}
          -  p_j^{1/2})^2 \bigg)^{1/2}\\ 
& = & \| Th^{1/2} \|_{2,Tr} \| h^{1/2}T - Th^{1/2} \|_{2,Tr}. 
\end{eqnarray*}
Now if $T$ happens to be a unitary operator then $\| Th^{1/2}
\|_{2,Tr} = \| h^{1/2} \|_{2,Tr} = 1$ and $\| h^{1/2}T - Th^{1/2}
\|_{2,Tr} = \| Th^{1/2}T^* - h^{1/2} \|_{2,Tr}$ and hence we can apply
the Powers-St{\o}rmer inequality after the inequalities above to get:
$$|Tr(h^{1/2}Th^{1/2}T^*) - tr(\phi(T)\phi(T^*))| \leq \| ThT^* -
h \|_{1,Tr}^{1/2}.$$ Finally, the Cauchy-Schwartz inequality
applied to the Hilbert-Schmidt operators implies that for every unitary 
operator $T \in B(H)$, 
\begin{eqnarray*}
tr(\phi(TT^*) - \phi(T)\phi(T^*)) & \leq & |1 - Tr(h^{1/2}Th^{1/2}T^*)| 
         + \| ThT^* - h \|_{1,Tr}^{1/2}\\ 
& = & |Tr(ThT^*) - Tr(h^{1/2}Th^{1/2}T^*)| + 
         \| ThT^* - h \|_{1,Tr}^{1/2}\\ 
& = & |Tr((Th^{1/2} - h^{1/2}T)h^{1/2}T^*)| + 
         \| ThT^* - h \|_{1,Tr}^{1/2}\\ 
& \leq & \| h^{1/2}T^* \|_{2,Tr} \| Th^{1/2} - h^{1/2}T \|_{2,Tr} 
         + \| ThT^* - h \|_{1,Tr}^{1/2}\\ 
& \leq & 2\| ThT^* - h \|_{1,Tr}^{1/2}.
\end{eqnarray*}
\end{proof}

\begin{thm}
\label{thm:invariantmeans}
Let $\tau$ be a tracial state on $A$.  Then the following are equivalent:
\begin{enumerate}
\item $\tau \in \TAim$.

\item There exists a sequence of u.c.p.\ maps $\phi_n : A \to M_{k(n)}
(\mathbb C)$ such that $\| \phi_n (ab) - \phi_n (a) \phi_n (b) \|_{2,tr}
\to 0$ and $\tau (a) = \lim_{n \to \infty} tr \circ \phi_n
(a)$, for all $a,b \in A$, where $\| x \|_{2,tr}^2 = tr(x^* x)$ for every
$x \in M_{k(n)} (\mathbb C)$.

\item For any faithful representation $A \subset B(H)$ there exists a
u.c.p.\ map $\Phi : B(H) \to \pi_{\tau}(A)^{\prime\prime}$ such that
$\Phi(a) = \pi_{\tau}(a)$.
\end{enumerate}
\end{thm}

\begin{proof}
(1) $\Longrightarrow$ (2). 
Let $A \subset B(H)$ be a faithful representation.  Since $\tau \in
\TAim$ we can find a state $\psi$ on $B(H)$ which extends $\tau$ and
such that $\psi(u T u^*) = \psi(T)$ for all unitaries $u \in A$ and
operators $T \in B(H)$. Since the normal states on $B(H)$ are dense in
the set of all states on $B(H)$ we can find a net of positive
operators $h_{\lambda} \in {\mathcal T}$ such that $Tr(h_{\lambda} T)
\to \psi(T)$ for all $T \in B(H)$.  Since $\psi(u^* T u) = \psi(T)$ it
follows that $Tr(h_{\lambda} T) - Tr((uh_{\lambda}u^*) T) \to 0$ for
every $T \in B(H)$ and unitary $u \in A$.  In other words,
$h_{\lambda} - uh_{\lambda}u^* \to 0$ in the weak topology on
${\mathcal T}$.  Hence, by the Hahn-Banach theorem, there are convex
combinations which tend to zero in L$^1$-norm.  In fact, taking
finite direct sums (i.e.\ considering n-tuples $(u_1 h_{\lambda} u_1^*
- h_{\lambda}, \ldots, u_n h_{\lambda} u_n^* - h_{\lambda})$) one
applies a similar argument to show that if ${\mathfrak F} \subset A$
is a finite set of unitaries then for every $\epsilon > 0$ we can find
a positive trace class operator $h \in {\mathcal T}$ such that $Tr(h)
= 1$, $|Tr(uh) - \tau(u)| < \epsilon$ and $\| h - uhu^* \|_1 <
\epsilon$ for all $u \in {\mathfrak F}$.  Since finite rank operators
are norm dense in ${\mathcal T}$ we may further assume that $h$ is
finite rank with rational eigenvalues. 

Applying Lemma \ref{thm:ozawalemma} to bigger and bigger finite sets
of unitaries and smaller and smaller epsilon's we can construct a
sequence of u.c.p.\ maps $\phi_n : B(H) \to M_{k(n)}({\mathbb C})$
such that tr($\phi_n (u)) \to \tau(u)$ and $|$tr$(\phi_n(uu^*) -
\phi_n(u)\phi_n(u^*))| \to 0$ for every unitary $u$ in a countable set
with dense linear span in $A$.  Since $\phi_n(uu^*) -
\phi_n(u)\phi_n(u^*) \geq 0$ we have $$ \| 1 - \phi_n(u)\phi_n(u^*) 
\|_{2,tr}^2 \leq \| 1 - \phi_n(u)\phi_n(u^*) \| tr(\phi_n(uu^*) -
\phi_n(u)\phi_n(u^*)) \to 0.$$  It follows that 
$\| \phi_n (ab) - \phi_n (a) \phi_n (b) \|_{2,tr}
\to 0$ for every $a, b \in A$.  Indeed, defining $\Phi = 
\oplus \phi_n : A \to \Pi M_{k(n)}({\mathbb C}) \subset 
l^{\infty}(R)$ we can compose with the natural quotient map
$l^{\infty}(R) \to R^{\omega}$ and it follows that every unitary such
that $\| \phi_n(uu^*) - \phi_n(u)\phi_n(u^*)
\|_{2,tr}^2 \to 0$ and $\| \phi_n(u^* u) - \phi_n(u^*)\phi_n(u)
\|_{2,tr}^2 \to 0$ will fall in the multiplicative domain of the 
composition.  However we have arranged that such unitaries have dense
linear span and hence all of $A$ falls in the multiplicative domain.

(2) $\Longrightarrow$ (3).  
Let $\phi_n : A \to M_{k(n)} (\mathbb C)$ be a sequence of u.c.p.\
maps with the properties stated in the theorem.  Identify each
$M_{k(n)} (\mathbb C)$ with a unital subfactor of $R$ and we can
define a u.c.p.\ map $A \to l^{\infty}(R)$ by $x \mapsto (\phi_n(x))
\in \Pi M_{k(n)}({\mathbb C}) \subset l^{\infty}(R)$.  Since the
$\phi_n$'s are asymptotically multiplicative in 2-norm one gets a
$\tau$-preserving $*$-homomorphism $A \to R^{\omega}$ by composing
with the quotient map $l^{\infty}(R) \to R^{\omega}$.  The weak
closure of $A$ under this mapping will be isomorphic to
$\pi_{\tau}(A)^{\prime\prime}$ and we can extend the mapping on $A$ to
all of $B(H)$ because (a) $l^{\infty}(R)$ is injective (hence we first
extend into $l^{\infty}(R)$) and (b) there exists a conditional
expectation $R^{\omega} \to \pi_{\tau}(A)^{\prime\prime}$.

(3) $\Longrightarrow$ (1).  Note that $A$ falls in the multiplicative
domain of $\Phi$ and hence $\Phi$ is a bimodule map; i.e.\ $\Phi(aTb)
= \pi_{\tau}(a)\Phi(T)\pi_{\tau}(b)$ for all $a,b \in A$ and $T \in
B(H)$.  From this observation one easily checks that if we let
$\tau^{\prime\prime}$ denote the vector trace on
$\pi_{\tau}(A)^{\prime\prime}$ then $\tau^{\prime\prime} \circ \Phi$
is a state on $B(H)$ which extends $\tau$ and which is invariant under
the action of the unitary group of $A$ on $B(H)$.  Hence $\tau$ is an
invariant mean. 
\end{proof}

\section{II$_1$-factor representations of $C^*({\mathbb F}_{\infty})$}

In this section we observe that every separable II$_1$-factor contains
a weakly dense copy of the universal C$^*$-algebra generated by a
countably infinite set of unitaries (i.e.\ $C^*({\mathbb
F}_{\infty})$).  Since every separable II$_1$-factor $M$ is generated
by a countable number of unitaries it follows from universality that
there is always a $*$-homomorphism $C^*({\mathbb F}_{\infty}) \to M$
with weakly dense range.  However, the next proposition completes the
II$_1$-factor representation theory of $C^*({\mathbb F}_{\infty})$; it
is not particularly deep but rather amounts to some universal
trickery.

\begin{prop}
\label{thm:mainlemma}
Let $M$ be a II$_1$-factor.  There exists a $*$-monomorphism $\rho :
C^*({\mathbb F}_{\infty}) \hookrightarrow M$ such that
$\rho(C^*({\mathbb F}_{\infty}))$ is weakly dense in $M$.
\end{prop}

\begin{proof} We first need to
write $C^*({\mathbb F}_{\infty})$ as an inductive limit of free
products of itself.  That is, we define $$A_1 = C^*({\mathbb
F}_{\infty}), A_2 = A_1 * C^*({\mathbb F}_{\infty}), \ldots, A_n =
A_{n-1} * C^*({\mathbb F}_{\infty}), \ldots,$$ where $*$ denotes the
full (i.e.\ universal) free product (with amalgamation over the scalars).  
Letting $A$ denote the inductive
limit of the sequence $A_1 \to A_2 \to \cdots$ it is easy to see (by
universal considerations) that $A \cong C^*({\mathbb F}_{\infty})$.
Since $A$ is residually finite dimensional (cf.\ \cite{davidson}) we
can find a sequence of integers $\{ k(n) \}$ and a unital
$*$-monomorphism $\sigma : A \hookrightarrow \Pi M_{k(n)} ({\mathbb
C})$.  Note that we may naturally identify each $A_i$ with a
subalgebra of $A$ and hence, restricting $\sigma$ to this copy of
$A_i$, get an injection of $A_i$ into $\Pi M_{k(n)} ({\mathbb C})$.

To construct the desired embedding of $A$ into $M$, it suffices to
prove the existence of a sequence of unital $*$-homomorphisms $\rho_i
: A_i \to M$ with the following properties:

\begin{enumerate}
\item Each $\rho_i$ is injective.

\item $\rho_{i+1}|_{A_i} = \rho_i$, where we identify $A_i$ with the
`left side' of $A_i * C^*({\mathbb F}_{\infty}) = A_{i+1}$.

\item The (increasing) union of $\{ \rho_i (A_i) \}$ is weakly dense
in $M$.
\end{enumerate}

To this end, we first choose an increasing sequence of projections
$p_1 \leq p_2 \leq \cdots $ from $M$ such that $\tau_M (p_i) \to 1$.
Then define orthogonal projections $q_2 = p_2 - p_1, q_3 = p_3 - p_2,
\ldots$ and consider the II$_1$-factors $Q_j = q_j M q_j$ for $j = 2,
3, \ldots$.  As is well known and not hard to construct, we can, for
each $j \geq 2$, find a unital embedding $\Pi M_{k(n)} ({\mathbb C})
\hookrightarrow Q_j \subset M$ and thus we get a sequence of
(orthogonal) embeddings $A \hookrightarrow \Pi M_{k(n)} ({\mathbb C})
\hookrightarrow Q_j \subset M$ which will be denoted by $\sigma_j$.

We are almost ready to construct the $\rho_i$'s. Indeed, for each $i
\in {\mathbb N}$ let $\pi_i : C^*({\mathbb F}_{\infty}) \to p_i M p_i$
be a (not necessarily injective!) $*$-homomorphism with weakly dense
range. We then define $\rho_1$ as $$\rho_1 = \pi_1 \oplus \bigg(
\bigoplus_{j \geq 2} \sigma_j|_{A_1} \bigg) : A_1 \hookrightarrow p_1
M p_1 \oplus \bigg( \Pi_{j \geq 2} Q_j \bigg) \subset M.$$ Note that
this is a unital $*$-monomorphism from $A_1$ into $M$ (since each
$\sigma_j$ is already faithful on all of $A$).  Now define a
$*$-homomorphism $\theta_2 : A_2 = A_1 * C^*({\mathbb F}_{\infty}) \to
p_2 M p_2$ as the free product of the $*$-homomorphisms $A_1 \to p_2 M
p_2$, $x \mapsto p_2 \rho_1 (x) p_2$, and $\pi_2 : C^*({\mathbb
F}_{\infty}) \to p_2 M p_2$.  We then put $$\rho_2 = \theta_2 \oplus
\bigg( \bigoplus_{j \geq 3} \sigma_j|_{A_2} \bigg) : A_2
\hookrightarrow p_2 M p_2 \oplus \bigg( \Pi_{j \geq 3} Q_j \bigg)
\subset M.$$ Note that $\rho_2|_{A_1} = \rho_1$.  Hopefully it is now
clear how to proceed.  In general, we construct a map (whose range is
dense in $p_{n+1} M p_{n+1}$) $\theta_{n+1}: A_{n+1} = A_n *
C^*({\mathbb F}_{\infty}) \to p_{n+1} M p_{n+1}$ as the free product
of the cutdown (by $p_{n+1}$) of $\rho_n$ and $\pi_{n+1}$.  This map
need not be faithful and hence we take a direct sum with $\oplus_{j
\geq n+2} \sigma_j|_{A_{n+1}}$ to remedy this deficiency.  It is then
easy to see that these maps have all the required properties and hence
the proof is complete.
\end{proof}

\section{Proof of main result} 

With Theorem \ref{thm:invariantmeans} and Proposition \ref{thm:mainlemma} in
hand we can now prove the main result.  

\begin{thm}
Let $M$ be a separable II$_1$-factor.  Then $M$ is embeddable into
$R^{\omega}$ if and only if $M$ has a weak expectation relative to
some weakly dense subalgebra.
\end{thm}

\begin{proof} ($\Longrightarrow$) First assume that 
$M \subset R^{\omega}$.  By Proposition \ref{thm:mainlemma} we may identify
$C^*({\mathbb F}_{\infty})$ with a weakly dense subalgebra of $M$.
Letting $\tau$ denote the unique trace on $M$ we first claim that
$\tau|_{C^*({\mathbb F}_{\infty})}$ is an invariant mean.  To see this
we note that since matrix algebras are weakly dense in $R$ we can find
a sequence $M_{k(n)}({\mathbb C}) \subset R$ such that each unitary in
$C^*({\mathbb F}_{\infty}) \subset M \subset R^{\omega}$ lifts to a
unitary in $\Pi M_{k(n)}({\mathbb C}) \subset l^{\infty}(R)$.  In
other words, there is a $*$-homomorphism $\sigma : C^*({\mathbb
F}_{\infty}) \to \Pi M_{k(n)}({\mathbb C})$ such that $\pi(\sigma(x))
= x$ for all $x \in C^*({\mathbb F}_{\infty})$, where $\pi :
l^{\infty}(R) \to l^{\infty}(R)/I_{\omega} = R^{\omega}$ is the
canonical quotient mapping.  By definition of the trace on
$R^{\omega}$ it follows that $\tau|_{C^*({\mathbb F}_{\infty})}$ is
the weak$-*$ limit of traces on matrix algebras composed with
homomorphisms $C^*({\mathbb F}_{\infty})
\to M_{k(n)}({\mathbb C})$ and hence $\tau|_{C^*({\mathbb
F}_{\infty})}$ is an invariant mean.  Now, if we move $M$ to its left
regular representation coming from $\tau$ then we can apply 
Theorem \ref{thm:invariantmeans} and
conclude that $M$ has a weak expectation relative to $C^*({\mathbb
F}_{\infty})$.

($\Longleftarrow$) Now suppose that there exists a weakly dense
C$^*$-algebra $A \subset M \subset B(H)$ and a u.c.p.\ map $\Phi :
B(H) \to M$ such that $\Phi(a) = a$ for all $a \in A$.  If $\tau$ is
the unique trace on $M$ then it follows that $\tau|_A$ is an invariant
mean just as in the proof of (3) $\Longrightarrow$ (1) from Theorem
\ref{thm:invariantmeans}.  From Theorem \ref{thm:invariantmeans}
 it follows that we can find a
sequence of u.c.p.\ maps $\phi_n : A \to M_{k(n)}({\mathbb C})$ which
are asymptotically multiplicative (in 2-norm) and which asymptotically
recover $\tau|_A$ after composing with the traces on
$M_{k(n)}({\mathbb C})$.  Hence the u.c.p.\ mapping $A \to
l^{\infty}(R)$ given by $x \mapsto (\phi_n(x)) \in \Pi
M_{k(n)}({\mathbb C}) \subset l^{\infty}(R)$ induces a
$\tau|_A$-preserving $*$-monomorphism $A \to R^{\omega}$ by composing
with the quotient map $l^{\infty}(R) \to R^{\omega}$.  It follows
(essentially due to uniqueness of GNS representations) that the weak
closure of $A$ in $R^{\omega}$ is isomorphic to $M$ and the proof is
complete.
\end{proof}

Finally we give the proof of the approximation property stated in the
introduction.  Note that a consequence of this result is that {\em if}
Connes' embedding problem is true (i.e.\ every separable II$_1$-factor
is embeddable) then $R$ is {\em the} basic building block for all
II$_1$-factors. It is hard for us to imagine that every II$_1$-factor
is built up from the inside by the nicest possible II$_1$-factor,
however a counterexample remains elusive.

\begin{cor}
Let $M \subset R^{\omega}$ be a II$_1$-factor, ${\mathfrak F} \subset
M$ be a finite set and $\epsilon > 0$ be given.  Then there exists a
complete order embedding $\Phi : R \hookrightarrow M$ (i.e.\ $\Phi$ is
an operator system isomorphism between $R$ and $\Phi(R)$ -- that is,
$\Phi$ is completely positive and $\Phi^{-1} : \Phi(R) \to R$ is also
completely positive) such that for each $x \in {\mathfrak F}$ there
exists $r \in R$ such that $\| x - \Phi(r) \|_2 < \epsilon$.
\end{cor}

\begin{proof} Let a finite set ${\mathfrak F} \subset
M$ and $\epsilon > 0$ be given.  Choose a projection $p \in M$ such
that $\tau(p) > 1 - \epsilon$.  Note that the corner $pMp$ is also
embeddable into $R^{\omega}$ (the fundamental group of $R^{\omega}$ is
${\mathbb R}_+$). 

Now let $C^*({\mathbb F}_{\infty}) \subset R$ be an identification
with a dense subalgebra of $R$ and $\pi : C^*({\mathbb F}_{\infty})
\hookrightarrow pMp$ be a $*$-monomorphism with weakly dense range.
By Theorem \ref{thm:invariantmeans} we can find a u.c.p.\ map 
$\Psi : R \to pMp$ which extends
$\pi$.  Since we can also find a unital $*$-homomorphism $\nu: R
\hookrightarrow (1-p)M(1-p)$ we get the desired complete order
embedding by defining $\Phi : R \to M$ by $\phi(r) = \Psi(r) \oplus
\nu (r)$.
\end{proof}

\bibliographystyle{amsplain}

\providecommand{\bysame}{\leavevmode\hbox to3em{\hrulefill}\thinspace}

\end{document}